\documentclass[final,1p,times]{elsarticle}

\usepackage{amssymb}
 \usepackage{amsthm}
\usepackage{amscd}
\usepackage{amsmath}
\usepackage{amsfonts}
\usepackage{amssymb}
\usepackage{graphicx}

\numberwithin{equation}{section}

\newcommand{\ra}{\rightarrow}
\newcommand{\p}{\partial}
\newcommand{\f}{\frac}

\newcommand{\be}{\begin{equation}}
\renewcommand{\ra}{\rightarrow}
\newcommand{\ee}{\end{equation}}
\newcommand{\bea}{\begin{eqnarray}}
\newcommand{\eea}{\end{eqnarray}}
\newcommand{\bna}{\begin{eqnarray*}}
\newcommand{\ena}{\end{eqnarray*}}

\renewcommand{\le}{\left}
\newcommand{\ri}{\right}

\journal{}

\begin{document}

\begin{frontmatter}

\title{Nonexistence of quasi-harmonic sphere with large energy}
\author[label1]{Jiayu Li}
 \ead[label1]{lijia@math.ac.cn}
\author[label2]{Yunyan Yang}
\ead[label2]{yunyanyang@ruc.edu.cn}
\address[label1]{Academy of Mathematics and Systems Science, Chinese Academy of Sciences, Beijing
100190, P. R. China}
\address[label2]{Department of Mathematics, Information School,
Renmin University of China, Beijing 100872, P. R. China}

\begin{abstract}
Nonexistence of quasi-harmonic spheres is necessary for long time
existence and convergence of harmonic map heat flows. Let $(N,h)$ be
a complete noncompact Riemannian manifolds. Assume the universal
covering of $(N,h)$ admits a nonnegative strictly convex function
with polynomial growth. Then there is no quasi-harmonic spheres
$u:\mathbb{R}^n\ra N$ such that
$$\lim_{r\ra\infty}r^ne^{-\f{r^2}{4}}\int_{|x|\leq r}e^{-\f{|x|^2}{4}}|\nabla u|^2dx=0.$$
This generalizes a result of the first named author and X. Zhu
(Calc. Var., 2009). Our method is essentially the Moser iteration
and thus very simple.

\end{abstract}


\end{frontmatter}

%\tableofcontents
\section{Introduction}
Let $(N,h)$ be a complete noncompact Riemannian manifolds. By the
Nash embedding theorem, there exists sufficiently large $K$ such
that $(N,h)$ is isometrically embedded in $\mathbb{R}^K$.  We say
that a map $w:\mathbb{R}^n\ra N\hookrightarrow\mathbb{R}^K$ is a
quasi-harmonic sphere if it satisfies
 \be\label{Quasi-equation}
 \Delta u=\f{1}{2}x\cdot\nabla u+A(u)(\nabla u,\nabla u),
 \ee
where $A(\cdot,\cdot)$ is the second fundamental form of $(N,h)$ in
$\mathbb{R}^K$. The Quasi-harmonic sphere arose from the study of
singularities of harmonic map heat flows \cite{Str,LW}. It is
closely related to the global smoothness and convergence of the
harmonic map heat flow \cite{CDW,DL,LiT}.

In \cite{LiZhu}, we have proved that if the universal covering of
$(N,h)$ admits a nonnegative strictly convex function with
polynomial growth, then there is no quasi-harmonic spheres of finite
energy, namely
 \be\label{finiteEnergy}\int_{\mathbb{R}^n}e^{-\f{|x|^2}{4}}|\nabla
u|^2dx<\infty.\ee The proof is based on the monotonicity inequality
for $u$ and John-Nirenberg inequality. In this note, we will use the
Moser iteration to prove a stronger result. Precisely we have the
following:
\\

\noindent{\bf Theorem 1.1} {\it Let $(N,h)$ be a complete Riemannian
manifold. Assume that $u$ is a quasi-harmonic sphere from
$\mathbb{R}^n$ $(n\geq 3)$ to $(N,h)$. Let
$(\widetilde{N},\widetilde{h})$ be the universal covering of
$(N,h)$. Suppose $(\widetilde{N},\widetilde{h})$ admits a
nonnegative strictly convex function $\widetilde{f}\in
C^2(\widetilde{N})$ with polynomial growth, i.e.
$\nabla^2\widetilde{f}(y)$ is positive definite for every $y\in
\widetilde{N}$ and $$\widetilde{f}(y)\leq
C(1+\widetilde{d}(y,y_0))^{2m}$$ for some $y_0\in\widetilde{N}$ and
positive integer $m$, where $\widetilde{d}(y,y_0)$ is the distance
between $y$ and $y_0$. If
\be\label{energy}\lim_{r\ra\infty}r^ne^{-\f{r^2}{4}}\int_{|x|\leq
r}e^{-\f{|x|^2}{4}}|\nabla u|^2dx=0,\ee then $u$ is a constant
map.}\\

We remark that if $u$ satisfies (\ref{energy}), then its energy
$\int_{\mathbb{R}^n}e^{-\f{|x|^2}{4}}|\nabla u|^2dx$ may be
infinite. In this sense, the conclusion of Theorem 1.1 is stronger
then that of \cite{LiZhu}. When $(N,h)$ is the standard real line
$\mathbb{R}$, the quasi-harmonic sphere becomes a quasi-harmonic
function, which is a solution to the equation
 $$
 \Delta u-\f{1}{2}x\cdot\nabla u=0 \quad{\rm in}\quad \mathbb{R}^n.
 $$
To prove Theorem 1.1, here we will use the Moser iteration instead
of using the monotonicity inequality for quasi-harmonic sphere and
the John-Nirenberg inequality for BMO space in \cite{LiZhu}.
Avoiding hard work from
harmonic analysis, our method looks very simple.\\

A special case of Theorem 1.1 is the following:\\

 \noindent{\bf Corollary 1.2} {\it Let $u$ be a quasi-harmonic
 function.
  If (\ref{energy}) is satisfied, then $u$ is a
 constant.}\\

 In view of Theorem 4.2 in \cite{LiWang}, any positive quasi-harmonic
 function $u:\mathbb{R}^n\ra \mathbb{R}$ with polynomial growth must be a constant. This is based
 on the gradient estimate. Its assumption can be interpreted by
 \be\label{energy2}\int_{|x|\leq r}e^{-\f{|x|^2}{4}}|\nabla u|^2dx\leq
 C(n)P(r),\ee
 where $C(n)$ is a universal constant and $P(r)$ is a polynomial
 with respect to $r$. Obviously the hypothesis (\ref{energy}) is
 much weaker than (\ref{energy2}). Hence the conclusion of Corollary
 1.2 is better than that of Theorem 4.2 in \cite{LiWang}.\\

 In the remaining part of this note, we will prove Theorem 1.1.

\section{Proof of Theorem 1.1}
Let $u:\mathbb{R}^n\ra (N,h)\hookrightarrow \mathbb{R}^K$ be a
quasi-harmonic sphere satisfying (\ref{Quasi-equation}). Denote
\be\label{wr}w(r)=\int_{\mathbb{S}^{n-1}}(2|u_r|^2-|\nabla
u|^2)d\theta=\int_{\mathbb{S}^{n-1}}(|u_r|^2-\f{1}{r^2}|u_\theta|^2)d\theta.\ee
It follows from (\ref{Quasi-equation}) that
 $\langle\Delta u,u_r\rangle=\f{r}{2}|u_r|^2$, and thus $\int_{\mathbb{S}^{n-1}}\langle\Delta u,u_r\rangle
 d\theta=\f{r}{2}\int_{\mathbb{S}^{n-1}}|u_r|^2d\theta$. Integration
 by parts implies
 \bea\label{basic-formula}
 \f{d}{dr}w(r)=\int_{\mathbb{S}^{n-1}}
 \le(\f{2}{r^3}|u_\theta|^2+\le(r-\f{2n-2}{r}\ri)|u_r|^2\ri)d\theta.
 \eea
 For details of deriving (\ref{basic-formula}), we refer the reader to
 \cite{LiZhu}.\\

 \noindent{\bf Lemma 2.1} {\it Let $w(r)$ be defined by (\ref{wr}), $w^+(r)$ be the positive part of $w(r)$,
 and $u$ be a quasi-harmonic sphere from $\mathbb{R}^n$ to $(N,h)$.
 Suppose
 \be\label{hypothesis}\int_0^re^{-\f{r^2}{4}}w^+(r)r^{n-1}dr\leq o(r^{-n}e^{\f{r^2}{4}})\quad{\rm as}\quad
 r\ra\infty.\ee
 Then there exists a constant $C$ depending only on $n$ and $w(2n)$
 such that
 $$\int_{\mathbb{B}_r}\le(d_N(u(x),u(0))\ri)^2dx\leq Cr^{n+1},$$
 where $d_N(\cdot,\cdot)$ denotes the distance function on
 $(N,h)$.}\\

 \noindent{\it Proof.} We can see from (\ref{basic-formula}) that $w^\prime(r)\geq
 0$ for $r\geq \sqrt{2n-2}$. We {\it claim} that $w(r)\leq 0$ for
 every $r\geq \sqrt{2n-2}$. Suppose not, there exists some $r_0\geq
 \sqrt{2n-2}$ such that $w(r_0)>0$. Then $w(r)\geq w(r_0)>0$ for every
 $r>r_0$ and
 \be\label{w'}
 w^\prime(r)\geq\le(r-\f{2n-2}{r}\ri)w(r).
 \ee
 We have by integrating $w^\prime(r)/w(r)$ from $r_0$ to $r$
  $$w(r)\geq w(r_0)r_0^{2n-2}e^{-\f{r_0^2}{2}}r^{2-2n}e^{\f{r^2}{2}}.$$
 Hence
 \bna\int_{r_0}^re^{-\f{r^2}{4}}w(r)r^{n-1}dr&\geq&w(r_0)r_0^{2n-2}e^{-\f{r_0^2}{2}}\int_{r_0}^r
 r^{1-n}e^{\f{r^2}{4}}dr\\
 &\geq&w(r_0)r_0^{2n-2}e^{-\f{r_0^2}{2}}r^{-n}\int_{r_0}^r
 re^{\f{r^2}{4}}dr\\
 &=&2w(r_0)r_0^{2n-2}e^{-\f{3r_0^2}{4}}r^{-n}e^{\f{r^2}{4}}.\ena
 This contradicts the assumption (\ref{hypothesis}) and thus confirms our claim.

 Now we estimate the growth order of the integral
 $\int_{\mathbb{B}_r}\le(d_N(u(x),u(0))\ri)^2dx$. For simplicity, we denote
 $d_N(u(x),u(0))$ by $d_N(x)$. In the polar coordinates in $\mathbb{R}^n$, we always identify $(r,\theta)$ with
 $x$. Notice that $d_N(r,\theta)\leq \int_0^r|u_r|ds$, one needs the
 following estimates, which can be obtained by using the H\"older inequality,
 the above claim, (\ref{wr}) and (\ref{w'}).
  \bna
 \int_{\mathbb{S}^{n-1}}\le(\int_0^r|u_r|ds\ri)^2d\theta&\leq&
 \int_{\mathbb{S}^{n-1}}r\le(\int_0^r|u_r|^2ds\ri)d\theta\\
 &\leq&r\int_0^{2n}\int_{\mathbb{S}^{n-1}}|u_r|^2d\theta ds+
 r\int_{2n}^{r}\int_{\mathbb{S}^{n-1}}|u_r|^2d\theta ds\\
 &\leq& Cr+r\int_{2n}^{r}\f{w^\prime(s)}{s-\f{2n-2}{s}}ds\\
 &\leq& Cr+Cr\int_{2n}^{r}w^\prime(s)ds\\
 &\leq& Cr+Cr(-w(2n))\\
 &\leq& Cr,
 \ena
 where $C$ is a constant depending only on $n$ and $w(2n)$.
 Hence we have
 \bna
  \int_{\mathbb{B}_r}d_N^2(x)dx&\leq&\int_0^rt^{n-1}\le\{
  \int_{\mathbb{S}^{n-1}}\le(\int_0^t|u_r|ds\ri)^2d\theta\ri\}
  dt\\
  &\leq& C\int_0^rt^ndt\leq Cr^{n+1}.
 \ena
 This concludes the lemma. $\hfill\Box$\\

 The following Lemma is elementary:\\

\noindent {\bf Lemma 2.2} {\it For every function $f$ defined on
$\mathbb{R}^n$,
 if there exists $k\in\mathbb{N}$ such that
 $$\int_{\mathbb{B}_r}|f(x)|dx\leq C_1r^k+C_2$$
 for some constants $C_1$ and $C_2$, then we have
 $$\int_{\mathbb{R}^n}e^{-\f{|x|^2}{4}} |f(x)|dx<\infty.$$}

 \noindent{\it Proof.} For sufficiently large $r$, it is easy to see
 that
 \bna
 \int_{\mathbb{R}^n\setminus\mathbb{B}_r}e^{-\f{|x|^2}{4}}
 |f(x)|dx&=&\sum_{j=1}^\infty\int_{\mathbb{B}_{2^jr}\setminus\mathbb{B}_{2^{j-1}r}}e^{-\f{|x|^2}{4}}
 |f(x)|dx\\
 &\leq&\sum_{j=1}^\infty e^{-4^{j-2}r^2}\int_{\mathbb{B}_{2^jr}}
 |f(x)|dx\\
 &\leq&\sum_{j=1}^\infty e^{-4^{j-2}r^2}(2^{kj}C_1r^k+C_2)\\
 &\leq&Cr^ke^{-\f{r^2}{4}}
 \ena
 for some constant $C$ depending only on $C_1$ and $C_2$. This immediately implies
 $$\lim_{r\ra\infty}\int_{\mathbb{R}^n\setminus\mathbb{B}_r}e^{-\f{|x|^2}{4}}
 |f(x)|dx=0,$$
 and thus gives the desired result. $\hfill\Box$\\

 We will use the Moser iteration of the following simple version (see for example Chapter 8 in \cite{GT}):\\

 \noindent{\bf Theorem A} {\it Let $u\geq 0$ be a weak solution of
 ${\rm div}(a\nabla u)\geq 0$ in $\mathbb{B}_{2\delta}(x_0)$, where $\delta>0$ is a constant, $x_0\in \mathbb{R}^n$,
 $a=a(x)$ satisfies
 $0<\lambda\leq a(x)\leq\Lambda$ in $\mathbb{B}_{2r}(x_0)$. Then for any $p>0$, there
 exists a constant $C$ depending only on $\Lambda/\lambda$, $n$ and $p$ such that
 $$\sup_{\mathbb{B}_{\delta}(x_0)}u
 \leq C\le(\f{1}{|\mathbb{B}_{2\delta}(x_0)|}\int_{\mathbb{B}_{2\delta}(x_0)}u^pdx\ri)^{1/p}.$$
 }

 For application of Theorem A, the following observation is crucial:\\

 \noindent{\bf Lemma 2.3} {\it Let $\rho(x)=e^{-|x|^2/4}$ on $\mathbb{R}^n$. Then for all $r>1$ and
 $x^\ast\in \overline{\mathbb{B}}_r=\{x\in \mathbb{R}^n: |x|\leq r\}$, there holds
 $$\sup_{x,\,y\,\in\mathbb{B}_{\f{2}{r}}(x^\ast)}\f{\rho(x)}{\rho(y)}\leq e^2.$$}

 \noindent{\it Proof.} Assume $x\in \mathbb{B}_{\f{2}{r}}(x^\ast)$. It
 is easy to see that
 $$\le(|x^\ast|-\f{2}{r}\ri)^2\leq |x|^2\leq \le(|x^\ast|+\f{2}{r}\ri)^2.$$
 Hence for $x,y\in \mathbb{B}_{\f{2}{r}}(x^\ast)$,
 \bna
 \f{\rho(x)}{\rho(y)}&\leq& \exp\le\{\f{1}{4}\le(|x^\ast|+\f{2}{r}\ri)^2-\f{1}{4}\le(|x^\ast|-
 \f{2}{r}\ri)^2\ri\}\\&\leq&\exp\le\{\f{2|x^\ast|}{r}\ri\}.
 \ena
 Note that $x^\ast\in \overline{\mathbb{B}}_r$, we get the desired
 result.$\hfill\Box$\\

 Now we are ready to prove Theorem 1.1 by using Theorem A.\\

 \noindent{\it Proof of Theorem 1.1.} Let $\widetilde{f}\in
C^2(\widetilde{N})$ be a nonnegative strictly convex function with
polynomial growth,
$u:\mathbb{R}^n\ra(N,h)\hookrightarrow\mathbb{R}^K$ be a
quasi-harmonic sphere, and $\widetilde{u}\in C^2(\widetilde{N})$ be
a lift of $u$.
 Define a function $\phi=\widetilde{f}\circ
 \widetilde{u}$. Let $\rho(x)=e^{-|x|^2/4}$. Then we have by a straightforward calculation
 \be\label{subelliptic}{\rm div}(\rho\nabla\phi)=\rho\nabla^2\widetilde{f}(\widetilde{u}(x))
 (\nabla \widetilde{u},\nabla \widetilde{u})\geq 0.\ee
 Assume $x^\ast\in \overline{\mathbb{B}}_r$ such that
 $\phi(x^\ast)=\sup_{\mathbb{B}_r}\phi$. It follows from the weak maximum principle for
 (\ref{subelliptic}) that
 $x^\ast\in\p\mathbb{B}_r$. By Lemma 2.3, we can apply Theorem A
 to the equation (\ref{subelliptic}) in the ball $\mathbb{B}_{\f{2}{r}}(x^\ast)$. This together with
 the hypothesis on $\widetilde{f}$ implies that for any $p>0$
 and $r>1$
 \bna
 \le(\f{1}{|\mathbb{B}_r|}\int_{\mathbb{B}_r}\phi^2dx\ri)^{1/2}&\leq&
 \sup_{\mathbb{B}_r}\phi\leq\sup_{\mathbb{B}_{\f{1}{r}}(x^\ast)}\phi\\
 &\leq&C\le(\f{1}{|\mathbb{B}_{\f{2}{r}}(x^\ast)|}\int_{\mathbb{B}_{\f{2}{r}}(x^\ast)}\phi^pdx\ri)^{1/p}\\
 &\leq&Cr^{n/p}\le(\int_{\mathbb{B}_{2r}}\phi^pdx\ri)^{1/p}\\
 &\leq&Cr^{n/p}\le(\int_{\mathbb{B}_{2r}}(1+\widetilde{d}\,^{2mp}(x))dx\ri)^{1/p},
 \ena
 where
 $\widetilde{d}(x)=\widetilde{d}_{\widetilde{N}}(\widetilde{u}(x),\widetilde{u}(0))$
 denotes the distance between $\widetilde{u}(x)$ and
 $\widetilde{u}(0)$ on the universal covering space $\widetilde{N}$ of
 $N$, $C$ is some constant depending only on $n$ and $p$. Clearly the assumption (\ref{energy}) implies (\ref{hypothesis}).
 Notice that Lemma 2.1 still holds when $u$ is replaced by
 $\widetilde{u}$, we have by choosing $p=1/m$ in the above inequality,
 $$\int_{\mathbb{B}_{r}}\phi^2dx\leq Cr^{n+(2n+1)2m},$$
 where $C$ is a constant depending only on $n$, $m$ and $\widetilde{u}$.
 From Lemma 2.2, we can see that
 \be\label{bound}\int_{\mathbb{R}^n}\rho\phi^2dx<\infty.\ee
 Take a cut-off function $\eta\in C_0^\infty(\mathbb{B}_{2r})$, $\eta\geq 0$ on $\mathbb{B}_{2r}$,
 $\eta\equiv
 1$ on $\mathbb{B}_{r}$, and $|\nabla \eta|\leq \f{4}{r}$. Testing
 the equation (\ref{subelliptic}) by $\eta^2\phi$, we obtain
 \bna
 \int_{\mathbb{R}^n}\eta^2\rho|\nabla
 \phi|^2dx&\leq&-\int_{\mathbb{R}^n}2\eta\phi\rho\nabla\eta\nabla\phi
 dx\\
 &\leq&2\le(\int_{\mathbb{R}^n}\eta^2\rho|\nabla
 \phi|^2dx\ri)^{1/2}\le(\int_{\mathbb{R}^n}\phi^2\rho|\nabla
 \eta|^2dx\ri)^{1/2}.
 \ena
 This together with (\ref{bound}) leads to
 $$\int_{\mathbb{B}_{r}}\rho|\nabla
 \phi|^2dx\leq \f{C}{r^2}$$
 for some constant $C$ depending only on the integral in (\ref{bound}).
 Passing to the limit $r\ra\infty$, we have $|\nabla \phi|\equiv 0$,
 which together with (\ref{subelliptic}) and that $\nabla^2\widetilde{f}$
 is positive definite implies that $|\nabla\widetilde{u}|\equiv 0$. Hence
 $\widetilde{u}$ is a constant
 map and thus $u$ is also a constant map. $\hfill\Box$\\

 {\bf Acknowledgements} The first author is partly supported by the National Science
 Foundation of China. The second author is partly supported by the NCET program.

\bigskip

\end{document}